\documentclass[12pt]{article}

\usepackage{amsmath,amsthm,amsfonts,latexsym,amscd}

\usepackage{eucal}

\usepackage{euler}
\usepackage{times}

\usepackage{xy}
\xyoption{all}

\swapnumbers

\theoremstyle{plain}
    \newtheorem{theorem}{Theorem}[section]
    \newtheorem{lemma}[theorem]{Lemma}
    \newtheorem{corollary}[theorem]{Corollary}
    \newtheorem{proposition}[theorem]{Proposition}

\theoremstyle{definition}
    
    \newtheorem{example}[theorem]{Example}
    \newtheorem{notation}[theorem]{Notation}

\theoremstyle{remark}
    \newtheorem{remark}[theorem]{Remark}

\numberwithin{equation}{section}

    \newcommand{\R}{\mathbb{R}}
    \newcommand{\C}{\mathbb{C}} 
    
    \newcommand{\Z}{\mathbb{Z}} 
    \newcommand{\Q}{\mathbb{Q}}

    \newcommand{\tensor}{\otimes}

    \DeclareMathOperator{\Kernel}{Kernel}
    \DeclareMathOperator{\Image}{Image}

\begin{document}

\title{Spaces with vanishing $\ell^2$-homology and their fundamental
groups (after Farber and Weinberger)}

\author{Nigel Higson, John Roe and Thomas Schick}

\maketitle

\section{Introduction}

Let $G$ be a finitely presented group.  Is $G$ the fundamental group
of some finite $CW$-complex whose integral homology is trivial
(meaning the homology is the same as that of a point)?  Elementary
homological algebra provides some simple  necessary conditions:
$$
        H_1(G, \mathbb Z) = 0 \quad \text{and} \quad H_2 (G,\Z )=0.
$$
An elegant observation of Kervaire~\cite{Kervaire} shows that these
necessary conditions are also sufficient:

\begin{theorem}
  Let $G$ be a finitely presented group and suppose that $H_1(G,
  \mathbb Z)$ and $ H_2 (G,\mathbb Z)$ are both zero. There is a
  connected $3$-dimensional finite $CW$-complex $X$ with $\pi_1(X)=G$
  such that $H_k(X,\mathbb Z)=0$ for all $k>0$.
\end{theorem}

Actually, Kervaire's attention was focussed on homology spheres, and
what he observed is this:

\begin{theorem}
  Let $G$ be a finitely presented group and suppose that $H_1(G,
  \mathbb Z)$ and $ H_2 (G,\mathbb Z)$ are both zero.  For every
  dimension $n\ge 5$ there is a homology $n$-sphere with fundamental
  group $G$.
\end{theorem}

The theorems rely on a computation, due to H.\ Hopf~\cite{Hopf}, of
the Hurewicz homomorphism in degree 2.  This and some other
ingredients have been used by Farber and Weinberger
\cite{Farber-Weinberger(1999)} to construct examples of finite
$CW$-complexes (3-dimensional) and smooth, closed manifolds
(6-dimensional) whose universal covers have vanishing
$\ell^2$-homology groups in all degrees.  See the survey article
\cite{Lott} for background material on $\ell^2$-homology and the `zero
in the spectrum question.'\footnote{The question, or conjecture, was
that no such examples of the sort that Farber and Weinberger
constructed could exist.}

The purpose of this note is to indicate that the same ingredients can
be used to prove direct analogues of Kervaire's theorems:

\begin{theorem}
\label{mainresult1} 
Let $G$ be a finitely presented group and suppose that the
homology groups $H_k(G,\ell^2(G))$ are zero for $k=0,1,2$. 
Then there is a  connected $3$-dimensional finite $CW$-complex $X$ with
$\pi_1(X)=G$ such that $H_k(X,\ell^2(G))=0$ for all
$k\ge 0$. 
\end{theorem}

\begin{theorem}
\label{mainresult2} 
Let $G$ be a finitely presented group and suppose that the
homology groups $H_k(G,\ell^2(G))$ are zero for $k=0,1,2$. 
For every dimension $n\ge 6$ there is a
closed manifold $M$ of dimension $n$ and with $\pi_1(M)=G$ such
that $H_k(M,\ell^2(G))=0$ for all $k> 0$.
\end{theorem}

\begin{remark} 
The difference between $n\ge 5$ in Kervaire's theorem and
$n\ge 6$ in ours is accounted for by the absence in the literature of
a suitable handle-cancellation lemma in $\ell^2$-homology.  Since our
main interest is in $CW$-complexes we shall not consider this issue
further here.
\end{remark}

The definitions of the $\ell^2$-homology groups $H_k(G,\ell^2(G))$ and
$H_k(X,\ell^2(G))$ will be reviewed in Section~\ref{homologicalsec}
below.  It is not hard to see that the hypotheses in the theorems can
not be reduced (compare \cite[11.3]{Lueck(1997)}), so that in fact the
vanishing of $H_k(G,\ell^2 (G))$ for $k=0,1,2$ is both necessary and
sufficient.

$C^*$-algebra theory makes an appearance in the Farber-Weinberger
proof, and the main novelty of our note is the more systematic use of
$C^*$-algebras, including homology with coefficients in a
$C^*$-algebra.  At other steps in the argument (notably the
construction of manifold examples from $CW$-examples) we have no
improvement to offer and we shall simply refer the reader to the paper
of Farber and Weinberger.

\section{Homological Preliminaries}
\label{homologicalsec}

Let $Z$ be a connected $CW$-complex and let $G$ be the fundamental
group of $Z$.  Form the cellular chain complex for the universal cover
of $Z$,
\begin{equation}
\label{complex}
        \xymatrix{ C _0(\widetilde Z) & \ar[l]_b C
        _1(\widetilde Z) & \ar[l]_b C _2(\widetilde
        Z)& \ar[l]_{\quad b} {\cdots}}, 
\end{equation}
which is a complex of projective (in fact free) right $\Z [G]$
modules.  The corresponding cellular homology groups $H_* (\widetilde
Z)$ are right $\Z[G]$-modules too, although not of course projective.

If $V$ is any left module over $\Z [G]$ then let us denote by $H_*(Z,
V)$ the homology of the tensor product complex
\begin{equation}
\label{tpcomplex}
        \xymatrix{ C _0(\widetilde Z)\otimes_{\mathbb Z[G]} V & \ar[l]_b C
        _1(\widetilde Z)\otimes_{\mathbb Z[G]} V & \ar[l]_b C _2(\widetilde
        Z)\otimes_{\mathbb Z[G]} V& \ar[l]_{\qquad\,  b} {\cdots}} .
\end{equation}

\begin{example} 
If $Z$ is a model for the classifying space $BG$ then $H_*(Z,V)$ is
isomorphic to the group homology of the module $V$:
$$
        H_*(BG, V)\cong H_*(G,V).
$$
Indeed, if $Z=BG$ then the cellular chain complex~(\ref{complex}) for
$\widetilde Z$ is a free resolution of the trivial $\Z[G]$-module
$\Z$.
\end{example}

There is an obvious \emph{coefficient homomorphism}
$$
\xymatrix{
h_n\colon H_n(\widetilde Z)\otimes_{\Z[G]} V \ar[r] & H_n(Z, V) }
$$
and we shall need the following simple observation concerning its
behaviour in degree 2:

\begin{lemma}
\label{hopflemma}
There is an exact sequence 
$$
\xymatrix{ H_2(\widetilde Z)\otimes_{\Z[G]} V \ar[r]^{\quad h_2} &
H_2(Z, V) \ar[r] & H_2(G, V)\ar[r]& 0 .  }
$$
In particular, if $H_2 (G, V)=0$ then the coefficient homomorphism is
surjective in degree 2.
\end{lemma}

\begin{proof}  
This is essentially the computation of Hopf to which we referred in
the introduction.  We shall prove the lemma in the case where $Z$ is
$2$-dimensional, which is the only case we shall need in this note.
The cellular chain complex for $Z$ may be prolonged in higher degrees
to obtain a projective resolution of the trivial $\Z[G]$-module $\Z$:
$$
\xymatrix{ 0& \Z\ar[l] &C _0(\widetilde Z)\ar[l] & \ar[l]_{b_1} C
        _1(\widetilde Z) & \ar[l]_{b_2} C _2(\widetilde
        Z)& \ar[l]_{\quad {b_3}} C_3 &  \ar[l]_{\quad {b_4}} C_4 &
        {\cdots}\ar[l]  }.
$$
Now the group $H_2 (Z,V)$ is equal to $\Kernel(b_2\otimes 1_V)$, while
the group $H_2 (\widetilde Z)\otimes _{\Z[G]} V$ is equal to
$\Kernel(b_2)\otimes _{\Z[G]} V$, which is equal to
$\Image(b_3)\otimes _{\Z[G]} V$.  Therefore the image of $H_2
(\widetilde Z)\otimes _{\Z[G]} V$ within $H_2 (Z,V)
=\Kernel(b_2\otimes 1_V) $ is $\Image(b_3\otimes 1_V)$.  The quotient
by this image is therefore
$$
        H_2(Z,V)/[H_2 (\widetilde Z)\otimes _{\Z[G]} V] =
        \Kernel(b_2\otimes 1_V)/\Image(b_3\otimes 1_V) = H_2 (G, V),
$$
as required.
\end{proof}

If $A$ is an auxiliary ring, and if $V$ is equipped with a right
$A$-module structure which commutes with the given left $\Z
[G]$-module structure then the complex~(\ref{tpcomplex}) is a complex
of right $A$-modules, and the homology groups $H_*(Z, V)$ have the
structure of right $A$-modules.  Moreover the coefficient homomorphism
is right $A$-linear.  We shall use this observation in the proof of
the main theorem, where we shall take $V$ to be a ring $A$ (considered
as a right module over itself), and where the $\Z[G]$-module structure
on $V$ will come from left multiplication by an embedded copy of $G$
in the group of invertible elements of $A$.

\section{Analytic Preliminaries}

We continue to denote by $Z$ a connected $CW$-complex with fundamental
group $G$.  The left regular representation of $G$ on the Hilbert
space $\ell^2 (G)$ provides $\ell^2(G)$ with the structure of a left
$\Z[G]$-module, and the homology groups $H_k (Z, \ell^2 (G))$ are the
$\ell^2$-homology groups to which we referred in the introduction.
They are denoted in a variety of different ways in the literature.

Notice that $\ell^2 (G)$ is a Banach space completion of the complex
group algebra $\C [G]$, or the real group algebra $\R [G]$ if we are
using real coefficients.\footnote{For the purposes of this note it
will be convenient to work over the reals.  We note that the real
cases of the results in this section follow immediately from their
complex counterparts by complexification.}  In the following section
it will be very convenient to work with completions which are not
merely Banach spaces but Banach \emph{algebras}.  For our purposes the
best choice is the \emph{reduced $C^*$-algebra of $G$}, denoted $C^*_r
(G)$, which is the norm-completion of $\C [G]$ in its left regular
representation as bounded operators on $\ell^2 (G)$. The group $G$
embeds into the group of invertible elements of $C^*_r (G)$, and left
multiplication by $G$ gives $C^*_r (G)$ the structure of a left
$\Z[G]$-module.  Hence we can form the groups $H_k (Z, C^*_r (G))$ as
in the previous section.  They are right $C^*_r (G)$-modules.

The purpose of this section is to prove the following result:

\begin{theorem}
\label{analytictheorem}
Suppose that $Z$ is a connected $CW$-complex with fundamental group
$G$ and finitely many cells in dimensions $0$ through $n$.  The
following are equivalent:
\begin{itemize}
\item The homology groups $H_k( Z, C^*_r(G))$ are zero in degrees $0$
through $n$.
\item The homology groups $H_k( Z, \ell^2(G) )$ are zero in degrees
$0$ through $n$.
\end{itemize}
\end{theorem}

The argument is an exercise in functional analysis, and to set the
proper context let us fix a unital $C^*$-algebra $A$ and a complex
$$\xymatrix{ E _0 & \ar[l]_b E _1& \ar[l]_b E _2& \ar[l]_b \cdots} $$
comprised of Hilbert $A$-modules and bounded, adjointable Hilbert
$A$-module \hbox{maps} (here and below, see Lance's text \cite{Lance(1995)},
especially Section~3, for information on Hilbert module theory).

\begin{lemma} 
The homology groups of the above complex vanish in degrees
$0$ through $n$ if and only if the `Laplace' operator $\Delta
=bb^*+b^*b$ is invertible on the spaces $E_0$ through $E_n$.
\end{lemma}

\begin{proof} 
Assume first $\Delta$ is surjective in degree $k$. For
arbitrary $x\in\Kernel( b)$ choose $y$ with $\Delta y=x$. Then $0=b\Delta y
= bbb^*y+bb^*by = bb^*by$. Therefore 
\begin{equation*} 
        0= \langle by,bb^*by\rangle = \langle b^*by,b^*by\rangle,
\end{equation*} 
and so $b^*by=0$. Consequently $x= bb^*y$, so that $x$ lies
in the image of $b$ and the $k$-th homology group vanishes.

If the homology groups in degree $k$ and $k-1$ vanish then the ranges
of $b_{k+1}\colon E_{k+1}\to E_k$ and $b_k\colon E_k\to E_{k-1}$ are
closed since they coincide with the kernels of the succeeding
differentials in our complex.  From \cite[Theorem 3.2]{Lance(1995)}
and its proof we obtain an orthogonal decomposition
$$
        E_k=\Image(b_{k}^*)\oplus \Image (b_{k+1})=\Kernel
        (b_k)^\perp \oplus \Kernel (b_{k+1}^*)^\perp .
$$ 
By the Open Mapping Theorem the operator $b_k$, and hence the operator
$b_k^*b_k^{\phantom{*}}$, is bounded below on the first summand, while
$b_{k+1}^{\phantom{*}}b_{k+1}^*$ is bounded below on the second (note
that since $\Image(b_{k+1})$ is closed, so is $\Image
(b_{k+1}^*)$). So the self-adjoint operator $\Delta$ is bounded below
on $E_k$ and is therefore invertible.
\end{proof}

Suppose now that $A$ is represented faithfully and non-degenerately on
a Hilbert space $H$ (we have in mind the regular representation of
$C^*_r (G) $ on $\ell^2 (G)$). Then we can form the Hilbert spaces
$E_k\otimes _A H$ by completing the algebraic tensor product over $A$
with respect to the norm associated to the  inner product
$$
        \langle e_1\otimes v_1 , e_2 \otimes v_2 \rangle_{E\otimes _A H} 
        = \langle v_1 , \langle e_1,e_2 \rangle_E v_2 \rangle_H.
$$
The Hilbert spaces so obtained assemble to form a complex of Hilbert
spaces and bounded linear maps.

\begin{lemma}  
Let $T$ be a bounded and adjointable operator on $E_k$.  The
operators $T$ on $E_k$ and $T\tensor I_H$ on $E_k\otimes_A H$ have the
same spectrum.  
\end{lemma}

\begin{proof} 
The map $T\mapsto T\otimes I$ is an injective homomorphism
of $C^*$-algebras and so preserves spectrum.  
\end{proof}

Putting the two lemmas together we obtain the following result:

\begin{lemma} 
\label{hilbertmodulelemma}
Suppose given a complex of Hilbert $A$-modules 
$$\xymatrix{ E _0 & \ar[l]_b E _1& \ar[l]_b E _2& \ar[l]_b \cdots} $$
and suppose that $H$ is a Hilbert space equipped with a faithful and
non-degenerate representation of $A$.  Then the following are
equivalent:
\begin{itemize}
\item The above complex has no homology in degrees $0$ through $n$.
\item The Hilbert module tensor product of the above complex with $H$
has no homology in degrees $0$ through $n$.\qed
\end{itemize}
\end{lemma}
 
Now if $A$ is a $C^*$-algebra with unit and if $E$ is a finitely
generated and projective module over $A$ (in the usual sense of
algebra) then $E$ may be given the structure of a Hilbert $A$-module,
and this structure is unique up to unitary isomorphism (this is a
generalization of the well-known fact that a complex vector bundle on
a compact space has an essentially unique Hermitian structure). See
\cite{Lance(1995)}.  Moreover all $A$-linear maps between such modules
are automatically bounded and adjointable.  Finally, $E$ is a finitely
generated projective $A$-module, and if $A$ is represented faithfully
and non-degenerately on $H$, then the algebraic and Hilbert module
tensor products $E\otimes _A H$ agree.  These observations allow us to
formulate a more algebraic version of the previous result:

\begin{lemma} 
Let $A$ be a unital $C^*$-algebra.  Suppose given a complex of
$A$-modules
$$
        \xymatrix{ E _0 & \ar[l]_{b_1} E _1& \ar[l]_{b_2} E _2&
        \ar[l]_{b_3} \cdots}
$$
(in the ordinary sense of algebra) for which $E_0$ through $E_n$ are
finitely generated and projective.  Suppose that $H$ is a Hilbert
space equipped with a faithful non-degenerate representation of $A$.
Then the following are equivalent:
\begin{itemize}
\item The above complex has no homology in degrees $0$ through $n$.
\item The \emph{algebraic} tensor product of the above complex over
$A$ with $H$ has no homology in degrees $0$ through $n$.
\end{itemize}
\end{lemma}

\begin{proof}  
Suppose first that the complex $E_*$ has zero homology in degrees $0$
through $n$.  The range of the differential $b_n\colon E_n \to
E_{n-1}$ is closed, since it is the kernel of $b_{n-1}$, and so by
\cite[Theorem 3.2]{Lance(1995)} the module $E_n$ splits as a direct
sum of $\Kernel (b_n)$ and its orthogonal complement.  In particular,
the module $\Kernel (b_n)$ is a direct summand of a finitely generated
module, and is therefore finitely generated itself.  It follows that
$\Image(b_{n+1})$ is finitely generated, and so we can find a finitely
generated free module $E_{n+1}'$ and a module map $f \colon E_{n+1}'
\to E_{n+1}$ for which the complex
$$
        \xymatrix{ E _0 & \ar[l]_{b_1} E _1& \ar[l]_{b_2} \dots &
        E_n \ar[l]_{b_{n}} & E_{n+1}' \ar[l]_{b_{n+1}f}}
$$
has vanishing homology in degrees $0$ through $n$.  Since all the
terms which appear now have Hilbert module structures we can appeal to
Lemma~\ref{hilbertmodulelemma} to conclude that the tensor product of
the displayed complex by $H$ has vanishing homology in degrees $0$
through $n$.  This implies the same vanishing result for the tensor
product of the original complex by $H$.

Suppose, conversely, that the algebraic tensor product complex
$E_*\otimes _A H$ has vanishing homology in degrees $0$ through $n$.
Then by the previous lemma, the homology of $E_*$ vanishes at least
though degree $n-1$, and so by the same argument as above, the module
$\Kernel (b_n)$ is still finitely generated.  We can therefore find
$f\colon E_{n+1}' \to E_{n+1}$ as in the first part of the proof, but
so that the tensor product complex involving $E'_{n+1}$ has vanishing
homology in degrees $0$ through $n$.  Once again,
Lemma~\ref{hilbertmodulelemma} now applies to complete the proof.
\end{proof}

\begin{proof}[Proof of Theorem~\ref{analytictheorem}.]
The algebraic tensor product of the module $C_k (\widetilde Z )
\otimes _{\Z[G]} C^*_r (G)$ over $C^*_r (G)$ with the Hilbert space
$\ell^2 (G)$ is $C_k (\widetilde Z ) \otimes _{\Z[G]} \ell^2 (G)$.  So
the theorem is an immediate consequence of the last lemma.
\end{proof}

Here are the consequences of the theorem that we shall need:

\begin{corollary} 
If $G$ is a finitely presented discrete group then the following are
equivalent: 
\begin{itemize}
\item  $H_k(G,\ell^2
(G))=0$ for $k=0,1,2$.
\item   $H_k(G,C^*_r 
(G))=0$ for $k=0,1,2$.\quad \qed 
\end{itemize}
\end{corollary}

\begin{remark} 
It is very easy to find examples of groups to which the proposition
above applies.  As Farber and Weinberger note, any threefold direct
product of finitely presented non-amenable groups will do the job.
\end{remark}

\begin{corollary} 
If $Z$ is a finite $CW$-complex with fundamental group $G$ then the
following are equivalent:
\begin{itemize}
\item  $H_k(Z,\ell^2
(G))=0$ for all $k\ge 0$.
\item   $H_k(Z,C^*_r 
(G))=0$ for all $k\ge 0$.\quad \qed 
\end{itemize}
\end{corollary}

\section{Proof of the Main Theorems} 

We shall prove the following result:

\begin{theorem}
\label{actualtheorem}
Let $G$ be a finitely presented group and suppose that the homology
groups $H_k(G,C^*_r(G))$ are zero for $k=0,1,2$.  Then there is a
connected, $3$-dimensional, finite $CW$-complex $X$ with $\pi_1(X)=G$
such that $H_k(X,C^*_r(G))=0$ for all $k\ge 0$.
\end{theorem}

In view of the results in the previous section, this is equivalent to
Theorem~\ref{mainresult1}.  Theorem~\ref{mainresult2} follows from
Theorem~\ref{mainresult1} by a regular neighborhood construction. We
refer the reader to the Farber-Weinberger paper
\cite{Farber-Weinberger(1999)} for details.

We shall work with the \emph{real} group $C^*$-algebra below, although
we note that the real and complex cases of Theorem~\ref{actualtheorem}
are easily derived from one another.

\begin{notation}
In the following, $G$ will denote a fixed finitely presented group for
which $H_k(G,C^*_r(G))=0$, for $k=0,1,2$.  All tensor products will be
taken over the ring $\Z [G]$ (in particular they will be algebraic ---
no Hilbert module tensor products will be involved).
\end{notation}

Nearly all of the argument below follows that of Farber and
Weinberger~\cite{Farber-Weinberger(1999)}, which in turn is inspired
by the argument of Kervaire mentioned earlier.

\begin{lemma}  
There exists a connected finite $2$-complex $Y$  with
fundamental group $G$ for which $H_2 (Y, C^*_r (G))$ is a finitely
generated and free $C^*_r (G)$-module.
\end{lemma}

\begin{proof}
Let $Z$ be a finite $2$-dimensional $CW$-complex with $\pi_1(Z)=G$
(for instance the presentation $2$-complex of our finitely presented
group $G$). Then in view of our assumptions on $G$, it follows that
$H_0(Z,C^*_r(G))$ and $ H_1(Z,C^*_r(G))$ are zero (compare the proof
of Lemma~\ref{hopflemma}).  Therefore, in view of the exact sequence
$$ 
\xymatrix@=6.5pt{ 0& \ar[l]C _0(\widetilde Z)\otimes C^*_r(G) &
   \ar[l]_b C _1(\widetilde Z)\otimes C^*_r(G) & \ar[l]_b C
   _2(\widetilde Z)\otimes C^*_r(G)& \ar[l]_b H_2(Z,C^*_r(G)) & \ar[l]
   0 ,}
$$
the module $H_2(Z,C^*_r(G))$ is finitely generated and stably free
(note that all the modules except perhaps the rightmost one are
actually free).  By wedging $Z$ with finitely many $2$-spheres we
obtain a finite $2$-complex $Y$ with fundamental group $G$ for which
$H_2(Y,C^*_r(G))$ is a free $C^*_r(G)$-module, as required.
\end{proof}

Define $j\colon H_2 (\widetilde Y)\to H_2 (\widetilde Y)\otimes C^*_r
(G)$ by associating to the class $x\in H_2 (\widetilde Y)$ the
elementary tensor $x\otimes 1\in H_2 (\widetilde Y)\otimes C^*_r (G)$.
Let us continue to denote by $h\colon H_2 (\widetilde Y)\otimes C^*_r
(G) \to H_2 ( Y, C^*_r (G))$  the coefficient map considered in
Section~\ref{homologicalsec}.

\begin{proposition}  
The image of the composition
$$
        h\circ j\colon H_2 (\widetilde Y) \to H_2(Y,C^*_r(G))
$$ 
contains a basis for the free $C^*_r(G)$-module $H_2(G,C^*_r(G))$.
\end{proposition}

\begin{proof} 
Since by Lemma~\ref{hopflemma} the coefficient map $h$ is surjective,
we can certainly find a basis $b_1,\dots b_d$ for $H_2(Y,C^*_r(G))$ in the
image of $h$.  Let us do so and write the basis elements as
$$
        b_k=\sum_{i=1}^{n_k} h(x_{i,k}\tensor a_{i,k}) =
        \sum_{i=1}^{n_k} h(x_{i,k}\tensor 1)a_{i,k}
$$ 
with $x_{i,k}\in H_2(\widetilde Y)$ and $a_{i,k}\in C^*_r(G)$ (in the
displayed formula we have used the fact that $h$ is
$C^*_r(G)$-linear).  Since $H_2(Y, C^*_r(G))$ is a free and finitely
generated module over $C^*_r(G)$ it has a Hilbert $C^*_r(G)$-module
structure.  The module multiplication operation
$$
        H_2(Y,C^*_r(G))\times C^*_r(G)\to H_2(Y,C^*_r(G))
$$ 
is continuous and the set of bases for $H_2(Y, C^*_r(G))$ is open
within the set of ordered $d$-tuples of elements in $H_2(Y,
C^*_r(G))$.  Therefore since $\Q [G]$ is dense in the real group
$C^*$-algebra $C^*_r(G)$, we can replace the elements $a_{i,k}\in
C^*_r(G)$ by sufficiently close $a_{i,k}'\in\Q [G]$ in such a way that
the elements
$$
        b_k' =\sum_{i=1}^{n_k}h(x_{i,k}\tensor 1)a_{i,k}'
$$ 
constitute a basis for $H_2(Y,C^*_r(G))$. We have now shown that the
map
$$
        H_2(\widetilde Y)\otimes \Q [G] \to H_2(Y,C^*_r(G))
$$
has within its image a basis for $H_2(Y,C^*_r(G))$.  But
multiplication with a non-zero real number is an automorphism of
$H_2(Y,C^*_r(G))$. So multiplying with the common denominator of all
coefficients of all $a_{i,k}'\in \Q [G]$ we obtain basis elements
$$
        b_k''=\sum_{i=1}^{n_k}h( x_{i,k}\tensor 1)a_{i,k}''
$$ 
with $a_{i,k}''\in \Z [G]$. Moreover, since
$$
        h(x_{i,k}\tensor 1)\!\cdot\!a_{i,k}''=h(x_{i,k}\tensor 1\!\cdot\!
        a_{i,k}'') =h(x_{i,k}\tensor a_{i,k}''\!\cdot\! 1) =
        h(x_{i,k}\!\cdot\! a_{i,k}''\otimes 1),
$$ 
the elements $b_k''$ lie in the image of $h\circ j$ as required.
\end{proof}

\medskip

\begin{proof}[Proof of Theorem~\ref{actualtheorem}.] 
Choose elements $v_1, \dots ,v_d \in H_2(\widetilde Y)$ which are sent
by $h\circ j$ to a basis for $H_2(Y,C^*_r(G))$. By the Hurewicz
Isomorphism Theorem each $v_k$ is represented by a map $S^2\to Y$.
Let us use these maps to attach $d$ $3$-cells to $Y$, and let us
denote by $X$ the $3$-dimensional $CW$-complex obtained in this
way. By construction,
$$
        H_0(X,C^*_r(G))=0 \quad \text{and}\quad H_1(X,C^*_r(G))=0.
$$ 
Moreover, the attaching maps are chosen exactly in such a way that the
differential $b\colon C_3(\widetilde X)\otimes C^*_r(G) \to
C_2(\widetilde X)\otimes C^*_r (G)$ gives an isomorphism between
$C_3(\widetilde X)\tensor C^*_r(G)$ and the space
$$ 
        H_2(Y,C^*_r(G))= \Kernel\bigl(b\tensor 1\colon C_2(\widetilde
        X)\tensor C^*_r(G) \to C_1(\widetilde X)\tensor C^*_r(G)\bigr).
$$ 
Thus we have killed all homology in degree 2 without creating any new
homology in degree 3.  Therefore we achieve
\begin{equation*}
  H_k(X,C^*_r(G))=0
\end{equation*}
for all $k$, as required.
\end{proof}


\end{document}